\documentclass[reqno,10pt,a4paper]{amsart}
\usepackage{srcltx}
\usepackage{graphicx,color,amssymb,latexsym,amsmath,amsfonts}
\usepackage{mathrsfs}

\numberwithin{equation}{section}
\theoremstyle{plain}
\newtheorem{keyword}{keyword}
\newtheorem{thm}{Theorem}[section]

\newtheorem{proposition}{Proposition}[section]
\newtheorem{ex}{Example}[section]
\begin{document}

\title{On One Property of Tikhonov\\
Regularization Algorithm}

\
\author{Mikhail Ermakov}\\

\maketitle

{erm2512@gmail.com}\\
{Institute of Problems of Mechanical Engineering, RAS and\\
St. Petersburg State University, St. Petersburg, RUSSIA\\
{Mechanical Engineering Problems Institute\\
Russian Academy of Sciences\\
Bolshoy pr.,V.O., 61\\
St.Petersburg\\
Russia}}
\vskip 0.3cm
St.Petersburg State University\\
University pr., 28, Petrodvoretz\\
198504 St.Petersburg\\
Russia
\vskip 0.3cm
St.Petersburg department of\\ Steklov
mathematical institute \\
Fontanka 27, St.Petersburg 191023

\begin{abstract} For linear inverse problem with Gaussian random noise we show that Tikhonov regularization algorithm is minimax in the class of linear estimators   and is asymptotically minimax in the sense of sharp asymptotic in the class of all estimators. The results are valid if some a priori information on a Fourier coefficients of solution is provided. For trigonometric basis this a priori information implies that the solution belongs to a ball in Besov space $B^r_{2\infty}$.
\end{abstract}

\begin{keyword}[class=AMS]
[Primary ]{65M30},{65R30},{62G08},{62J07}
\end{keyword}

\begin{keyword}
Tikhonov regularization algorithm, linear inverse problem, nonparametric estimation, asymptotic minimaxity, asymptotic efficiency
\end{keyword}



\section{\bf  Introduction \label{s1}}
Tikhonov regularization algorithm (TRA) is very popular \cite{tik, ph} thanks to many remarkable properties. We mention only two of them.  TRA is minimax for deterministic noise \cite{mik} and is Bayes estimator \cite{tia, en} in the problems with Gaussian random noise and Gaussian a priori measure.
In paper  we explore the minimax properties of TRA in linear inverse problems with  Gaussian random noise.  We show that TRA is minimax in the class of linear estimator and asymptotically minimax in the class of all estimators. In these setups a priori information is provided that Fourier coefficients of solution  satisfies the same  restrictions as Fourier coefficients of functions in a ball in Besov space  $B^s_{2\infty}$ for the case of trigonometric basis.

Such a form of a priori information is rather natural.
\vskip 0.3cm
This is rather reasonable information on a solution smoothness.
\vskip 0.3cm
For the most nonparametric estimators these sets are the largest sets with a given rates of convergence \cite{ker93, ker02}.
\vskip 0.3cm
For linear statistical estimators these sets are the largest sets  with a given rate of convergence \cite{rio}.
\vskip 0.3cm
The asymptotic minimaxity of TRA is proved in the sense of sharp asymptotic. The asymptotically minimax nonparametric estimators in the sense of sharp asymptotic has been obtained earlier, only if a priori information is provided that a solution  belongs to ellipsoid in $L_2$, in particular, a ball in Sobolev space \cite{pin, ts, jo}.

There are numerous research on sharp adaptive minimax estimation \cite{en, jo, ca,ts, mair}.
The results on adaptive estimation in Pinsker model \cite{jo,ts} are easily carried over on  paper setup.

In what follows we shall denote letters $c, C$ positive constants and $ a_\epsilon \asymp b_\epsilon$ implies $c < a_\epsilon/b_\epsilon < C$.
\section{\bf Main Results}
Let $H$ be separable Hilbert space and let $A: H \to H$ be known self-adjoint linear bounded operator.

We wish to estimate a solution of linear equation
$$
f = Ax, \qquad x \in H,
$$
on observation $ Y = f + \epsilon\xi$ where $\xi$  is Gaussian random error and $\epsilon > 0$  defines the level of noise.

Let $\{a_i\}_{i=1}^\infty$ and $\{\phi_i\}_{i=1}^\infty$ be eigenvalues and eigenvectors  of operator $A$ respectively. Then we can rewrite the vector $Y = \{y_i\}_{i=1}^\infty$ in the
following form

\begin{equation}\label{r1}
y_i = a_i x_i + \epsilon \sigma_i \xi_i
\end{equation}
where $x_i= <x,\phi_i>$, $y_i = <Y,\phi_i>$ and $ \xi_i = \sigma_i^{-1}<\xi,\phi_i>, 1 \le i < \infty$, $Var[\xi_i] =1$. We suppose that $\xi_i, 1 \le i < \infty$ are i.i.d.r.v.'s and $E[\xi_i] =0$.
The representation (\ref{r1}) holds in particular if $\xi$ is Gaussian white noise.

Here $<a,b>$ denote inner product of vectors $a,b \in H$. For any $a  \in H$ denote $||a|| = <a,a>^{1/2}$.

Suppose a priori information is provided that
\begin{equation}\label{r2}
x \in B^r_{2\infty} = \left\{\theta= \{\theta_i\}_{i=1}^\infty : \sup_k k^{2r} \sum_{i=1}^\infty \theta_i^2 \le P_0\right\}
\end{equation}
with $r>0$.

We say that linear estimator $\hat x_\epsilon = \{\hat x_{\epsilon j}\}_{j=1}^\infty$ is minimax linear estimator if
\begin{equation}\label{r3}
\sup_{x \in B^r_{2\infty}} E_x||\hat x_\epsilon - x ||^2 =
\inf_{\lambda}\sup_{x \in B^r_{2\infty}} E_x||\hat x_{\epsilon\lambda} - x ||^2.
\end{equation}
where $\lambda= \{\lambda_j\}_{j=1}^\infty$, $\hat x_{\epsilon\lambda} = \{\hat x_{\epsilon j\lambda_j}\}_{j=1}^\infty, \hat x_{\epsilon j\lambda_j} = \lambda_j y_j, \lambda_j \in R^1, 1 \le j < \infty.$

We say that the estimator $\hat x_\epsilon$ is asymptotically minimax if
\begin{equation}\label{}
\sup_{x \in B^r_{2\infty}} E_x||\hat x_\epsilon - x ||^2 =
\inf_{\tilde x_\epsilon \in \Psi}\sup_{x \in B^r_{2\infty}} E_x||\tilde x_\epsilon - x ||^2(1 + o(1))
\end{equation}
as $\epsilon \to 0$. Here $\Psi$ is the set of all estimators.

The minimaxity of TRA in the class of linear estimators will be proved if the following assumption holds.

\noindent{\bf A}. For all $j>1$
\begin{equation}\label{}
\frac{\sigma_j^2 a_{j-1}^2((j-1)^{-2r}- j^{-2r})}{\sigma_{j-1}^2 a_j^2(j^{-2r}- (j+1)^{-2r})} > 1.
\end{equation}
\begin{thm}\label{t1} Assume A. Then TRA is  minimax on the set of all linear estimators with
\begin{equation}\label{}
\lambda_j = a_j^{-1}\frac{a_j^2P_0(j^{-2r}- (j+1)^{-2r})}{a_j^2P_0(j^{-2r}- (j+1)^{-2r}) + \epsilon^2\sigma_j^2}.
\end{equation}
The asymptotically minimax risk equals
\begin{equation}\label{}
R_{l\epsilon} = \epsilon^2(1+o(1))\sum_{j=1}^\infty \frac{\sigma_j^2P_0(j^{-2r}- (j-1)^{-2r})}{a_j^2P_0(j^{-2r}- (j-1)^{-2r}) + \epsilon^2\sigma_j^2 }.
\end{equation}
\end{thm}
The asymptotic minimaxity of TRA will be proved if the following assumptions hold.

\noindent{\bf B1} For $j>j_0$, there holds $|a_j/a_{j+1}| \ge 1$.

\noindent{\bf B2} There holds $0 < c< \sigma_j^2 < C < \infty$.

\noindent{\bf B3}. For all $j>j_0$
\begin{equation}\label{}
\frac{\sigma_j^2 a_{j-1}^2j^{2r+1}}{\sigma_{j-1}^2 a_j^2(j-1)^{2r+1}} > 1.
\end{equation}
\begin{thm}\label{t2} Assume B1-B3. Then TRA is asymptotically minimax on the set of all estimators with
\begin{equation}\label{}
\lambda_j = a_j^{-1}\frac{a_j^2}{a_j^2 + (2rP_0)^{-1}\epsilon^2\sigma_j^2j^{2r+1}}.
\end{equation}
The asymptotically minimax risk equals
\begin{equation}\label{}
R_\epsilon = \epsilon^2(1+o(1))\sum_{j=1}^\infty \frac{\sigma_j^2}{a_j^2 + (2rP_0)^{-1}\epsilon^2\sigma_j^2 j^{2r+1}}.
\end{equation}
\end{thm}
\begin{ex}\label{e1} Let $|a_j| = C j^{-\gamma}(1 + o(1))$ and $\sigma_j =1$. Then
\begin{equation}\label{}
R_\epsilon= \epsilon^{\frac{4r}{1+2r+2\gamma}}\frac{\pi}{2r
\sin\left(\frac{\pi(2\gamma+1)}{2r}\right)}(2rP_0)^{\frac{2\gamma+1}
{2\gamma+2r+1}}C^{-\frac{2r}{2\gamma+2r+1}}(1+o(1)).
\end{equation}
\end{ex}
\begin{ex}\label{e1} Let $|a_j| = Cj^{-\alpha}\exp\{-Bj^\gamma\}$ and $\sigma_j =1$. Then
\begin{equation}\label{}
R_\epsilon= |\log \epsilon|^{-2r/\gamma}P_0 B^{2r/\gamma}(1 + o(1)).
\end{equation}
\end{ex}
\section{\bf Proof of Theorem \ref{t1}}
We begin with the proof of lower bound. Denote $\theta_j^2 = P_0(j^{-2r}- (j+1)^{-2r}), \theta= \{\theta_j\}_{j=1}^\infty$.

We have
\begin{equation}\label{}
\inf_\lambda\sup_{x\in B^r_{2\infty}} E_x || \hat x_\lambda- x||^2 \ge \inf_\lambda E_\theta ||\hat\theta_\lambda - \theta||^2\\=
\epsilon^2\sum_{j=1}^\infty\frac{\theta_j^2\sigma_j^2}{\theta_j^2a_j^2 + \epsilon^2\sigma_j^2}
\end{equation}
and infimum is attained for
$$
\lambda_j = a_j^{-1}\frac{a_j^2\theta_j^2}{\theta_j^2a_j^2 + \epsilon^2\sigma_j^2}.
$$
Proof of upper bound is based on the following reasoning. Let $x=\{x_j\}_{j=1}^\infty \in B^r_{2\infty}$. For all $k$ denote
$$ u_k = k^{2r}\sum_{j=k}^\infty x_j^2.$$
Then $x_k^2= k^{-2r} u_k - (k+1)^{-2r} u_{k+1}$.

For the sequence of $\lambda_j$ defined in Theorem \ref{t1}, we have
\begin{equation}\label{ux1}
\begin{split}&
E_x\sum_{j=1}^\infty (\lambda_jy_j - x_j)^2 =\epsilon^2 \sum_{j=1}^\infty \lambda_j^2\sigma_j^2 a_j^{-2}+ \sum_{j=1}^\infty(1 -a_j\lambda_j)^2x_j^2\\&= \epsilon^2 \sum_{j=1}^\infty \lambda_j^2\sigma_j^2 a_j^{-2} +  \sum_{j=1}^\infty\left(\frac{1}{\theta_j^2\sigma_j^{-2}a_j^2\epsilon^{-2} +1}\right)^2(j^{-2r} u_j - (j+1)^{-2r} u_{j+1})\\&=\epsilon^2 \sum_{j=1}^\infty \lambda_j^2\sigma_j^2 a_j^{-2} + \left(\frac{1}{\theta_1^2\sigma_1^{-2}a_1^2\epsilon^{-2} +1}\right)^2 u_1
\\&-
\sum_{j=2}^\infty u_j j^{-2r}\left((\theta_{j-1}^2\sigma_{j-1}^{-2}a_{j-1}^2\epsilon^{-2} +1)^{-2} -
(\theta_j^2\sigma_j^{-2}a_j^2\epsilon^{-2} +1)^{-2}\right).
\end{split}
\end{equation}
By A, the last addendums in the right hand-side of (\ref{ux1}) are negative. Therefore the supremum of right hand-side of (\ref{ux1}) is attained for $u_j= P_0$, $1 \le j< \infty$. This completes the proof of Theorem  (\ref{t1}).
\section{\bf Proof of Theorem \ref{t2}}
The upper bound follows from Theorem \ref{t1}. Below the proof of lower bound will be provided. This proof has a lot of common feachers with the proof of lower bound in Pinsker Theorem  \cite{jo, pin, ts}.

Fix values $\delta, 0<\delta<1,$ and $\delta_1, 0 < \delta_1< P_0$. Define a family of natural  numbers $k_\epsilon, \epsilon>0,$ such that $a_{k_\epsilon}^2\epsilon^{-2} \sigma^2_{k_\epsilon} P_0(2r)^{-1}k_\epsilon^{-2r-1} = 1 + o(1)$ as $\epsilon \to 0$. Define sequence $\eta = \{\eta_j\}_{j=1}^\infty$ of Gaussian i.i.d.r.v.'s $\eta_j = \eta_{j\delta\delta_1}, E[\eta_{j}] = 0, \mbox{Var}[\eta_j] = (P_0 - \delta_1)(2r)^{-1} j^{-2r-1}$, if $\delta k_\epsilon \le j \le \delta^{-1} k_\epsilon$, and $\eta_j=0$ for $j < \delta k_\epsilon$ and $j > \delta^{-1}k_\epsilon$

Denote $\mu$ the probability measure of random vector $\eta$. Define $\tilde x$ Bayes estimator with a priory measure $\mu$.

 Define the conditional probability measure $\nu_\delta$ of random vector
  $\eta$ given $\eta \in  B^r_{2\infty}(P_0).$
  Define $\bar x $ Bayes estimator $x$ with a priori measure $\nu_\delta$.
 Denote $\theta$ the random variable having probability measure $\nu_\delta$.

For any estimator $\hat x$ we have
\begin{equation}\label{e1}
\begin{split}&
\sup_{x \in B^s_{2\infty}} E_x||\hat x -x||^2 \ge E_{\nu_\delta}E_\theta ||\hat x - \theta||^2\\&\ge
E_\mu E_\eta ||\tilde x -\eta||^2 -E_\mu E_\eta(||\bar x - \eta||^2, \eta \notin B^r_{2\infty}) P^{-1}_\mu(\eta \in B^r_{2\infty}) .
\end{split}
\end{equation}
We have
\begin{equation}\label{e2}
E_\mu E_\eta ||\tilde x -\eta||^2 = \epsilon^2(1+o(1))\sum_{j=l_1}^{l_2}
\frac{\sigma_j^2}{a_j^2 + (2r(P_0-\delta_1))^{-1}\epsilon^2\sigma_j^2 j^{2r+1}}\doteq I(P_0-\delta_1)
\end{equation} with $l_1 = [\delta k_\epsilon]$ and $l_2 =[\delta^{-1}k_\epsilon].$ Here $[a]$ denotes whole part of a number $a \in R^1$.

Since $$||\bar x ||^2 \le \sup_{x\in B^r_{2\infty}} ||x||^2 \le P_0,$$ we have
\begin{equation}\label{e3}
\begin{split}&
E_\mu E_\eta(||\bar x - \eta||^2, \eta \notin B^r_{2\infty}) \le  2 E_\mu E_\eta ( ||\bar x||^2 + ||\eta||^2, \eta \notin B^r_{2\infty}) \\& \le 2 P_0 P_\mu(\eta \notin B^r_{2\infty}) + \sum_{j=l_1}^{l_2} (E_\mu \eta_j^4)^{1/2} P_\mu^{1/2}(\eta \notin B^r_{2\infty}).
\end{split}
\end{equation}
Since
$E_\mu [\eta_j^4] \le C j^{-2r-2},$ we have
\begin{equation}\label{e4}
\sum_{j=l_1}^{l_2} (E_\mu \eta_j^4)^{1/2} \le C\delta^{-r} k_\epsilon^{-2r}.
\end{equation}
It remains to estimate
\begin{equation}\label{e5}
P_\mu(\eta \notin B^r_{2\infty}) = P(\max_{l_1 \le i \le l_2} i^{2r} \sum_{j = i}^{l_2} \eta_j^2-P_0(1-\delta_1/2) > P_0\delta_1/2) \le \sum_{i=l_1}^{l_2} J_i
\end{equation}
with
$$
J_i = P\left( i^{2r} \sum_{j = i}^{l_2} \eta_j^2-P_0(1-\delta_1/2)> P_0\delta_1/2\right)
$$
To estimate $J_i$ we implement the following Proposition \cite{hs}
\begin{proposition}\label{p1} Let $\xi = \{\xi_i\}_{i=1}^l$ be Gaussian random vector with i.i.d.r.v.'s $\xi_i$, $E[\xi_i] = 0, E[\xi_i^2]=1$. Let $A\in R^l\times R^l$ and $\Sigma = A^T A$. Then
\begin{equation}\label{e6}
P(||A\xi||^2 > \mbox{tr}(\Sigma) + 2\sqrt{\mbox{tr}(\Sigma^2)t} + 2 ||\Sigma||t) \le \exp\{-t\}.
\end{equation}
\end{proposition}
We put $\Sigma= \{\sigma_{lj}\}_{l,j=i}^{l_2}$ with $\sigma_{jj} = j^{-2r-1}i^{2r}\frac{P_0-\delta}{2r}$ and $\sigma_{lj} =0$ if $l\ne j$. Then
\begin{equation}\label{e7}
2\sqrt{\mbox{tr}(\Sigma^2)t} + 2 ||\Sigma||t =  \frac{P_0-\delta}{r(4r+1)} \sqrt{i^{-1}t}(1+o(1)) + i^{-1}t \doteq V_i(t)
\end{equation}
Since $i > \delta k_\epsilon$, we can put $t = k_\epsilon^{1/2}$. Then $V_i(t) < Ck_\epsilon^{-1/2}$ and implementing (\ref{e6})  we have
\begin{equation}\label{e8}
J_i < \exp\{-k_\epsilon^{-1/2}\}
\end{equation}
and therefore
\begin{equation}\label{e8}
\sum_{j=l_1}^{l_2} J_i \le \delta^{-1}k_\epsilon\exp\{- k_\epsilon^{1/2}\}
\end{equation}
To complete the proof it remains to estimate $R_\epsilon-I(P_0-\delta_1)$.

By straightforward estimation, it is easy to verify that
\begin{equation}\label{e9}
|I(P_0) - I(P_0 -\delta_1)| < C\delta_1 I(P_0)
\end{equation}
We have
\begin{equation}\label{e10}
\begin{split}&
\epsilon^2\sum_{j=1}^{l_1} \frac{\sigma_j^2}{a_j^2 + (2rP_0)^{-1}\epsilon^2\sigma_j^2 j^{2r+1}} \asymp \epsilon^2\sum_{j=1}^{l_1} \sigma_j^2 a_j^{-2} \\& <C\delta_1\sum_{j=l_1}^{k_\epsilon}\sigma_j^2 a_j^{-2} \asymp C\delta_1 \epsilon^2\sum_{j=l_1}^{k_\epsilon} \frac{\sigma_j^2}{a_j^2 + (2rP_0)^{-1}\epsilon^2\sigma_j^2 j^{2r+1}}
\end{split}
\end{equation}
We have
\begin{equation}\label{e11}
\begin{split}&
\epsilon^2(1+o(1))\sum_{j=l_2}^\infty \frac{\sigma_j^2}{a_j^2 + (2rP_0)^{-1}\epsilon^2\sigma_j^2 j^{2r+1}} \asymp \sum_{j=l_2}^\infty j^{-2r-1}\\& \le
\delta^{2r}C  \sum_{k_\epsilon}^{l_2} j^{-2r-1} \asymp \delta^{2r}C  \sum_{k_\epsilon}^{l_2} \frac{\sigma_j^2}{a_j^2 + (2rP_0)^{-1}\epsilon^2\sigma_j^2 j^{2r+1}}.
\end{split}
\end{equation}

\end{document}